\documentclass{ifacconf}

\usepackage{graphicx, bm}      
\usepackage{subcaption}
\usepackage{natbib}        

\begin{document}
\begin{frontmatter}

\title{Boundary Control for Suppressing Chaotic Response to Dynamic Hydrogen Blending in a Gas Pipeline}

\thanks[footnoteinfo]{This study was funded by the LANL Laboratory Directed R\&D Project ``Efficient Multi-scale Modeling of Clean Hydrogen Blending in Large Natural Gas Pipelines to Reduce Carbon Emissions'', the U.S. Department of Energy through the LANL/LDRD Program, and the Center for Non-Linear Studies. Research conducted at Los Alamos National Laboratory is done under the auspices of the National Nuclear Security Administration of the U.S. Department of Energy under Contract No. 89233218CNA000001. Report No. LA-UR-23-32580.}

\author[First]{Luke S. Baker} 
\author[Second]{\& Anatoly Zlotnik} 

\address[First]{Center for Nonlinear Studies, Theoretical Division,
    Los Alamos National Laboratory, 
    Los Alamos, NM USA (e-mail: lsbaker@lanl.gov)}
\address[Second]{Applied Mathematics \& Plasma Physics Group,
    Los Alamos National Laboratory, 
    Los Alamos, NM  USA (e-mail: azlotnik@lanl.gov)}

\begin{abstract}                
It is known that periodic forcing of nonlinear flows can result in a chaotic response under certain conditions.  Such non-periodic and chaotic solutions have been observed in simulations of heterogeneous gas flow in a pipeline with periodic, time-varying boundary conditions. In this paper, we examine a proportional feedback law for boundary control of a parabolic partial differential equation system that represents the flow of two gases through a pipe.  We demonstrate that periodic variation of the mass fraction of the lighter gas at the pipe inlet can result in the chaotic propagation of gas pressure waves, and show that  appropriate flow control can suppress this response.    We examine phase space solutions for the single pipe system subject to boundary control, and use numerical experiments to characterize conditions for the controller gain to suppress chaos.
\end{abstract}


\end{frontmatter}

\section{Introduction}

Chaotic behavior has been observed in a variety of continuum systems that involve gas dynamics, which are defined by partial differential equation (PDE) systems.  Prominent examples in nature include seasonal resonances and weather patterns in the oceanic atmosphere related to cyclical wind bursts ~(\cite{tziperman1994nino, eisenman2005westerly}). Examples in engineered systems include flame combustion of mixtures of hydrogen and air in microchannels ~(\cite{pizza2008dynamics, alipoor2016combustion}), along with the combustion process of a premixed natural gas engine ~(\cite{ding2017analysis}).  The chaotic dynamics in these systems are observed as high sensitivity to small perturbations in initial conditions that result in rapid divergence of trajectories, and as irregular responses to periodic boundary conditions.  

Although chaos can be directly exploited in certain applications, such as cryptography ~(\cite{kocarev2001chaos}), it is in general not desirable in most engineering settings, and in particular those that involve fluid flows.  The suppression of chaotic dynamics is desired in engineered systems with observable states and controllable inputs, in order to maintain predictable function within design limits.  The application of control methods to suppress chaos in a dynamical system has been examined in various studies ~(\cite{boccaletti2000control, zhang2009controlling}).  Approaches including open-loop control ~(\cite{dudnik1983synchronization}), feedback control ~(\cite{pyragas1992continuous, hu1995feedback, schuster1997control}), and the Ott-Grebogi-Yorke (OGY) method ~(\cite{ott1990controlling}) have all demonstrated the ability to control or suppress chaotic behavior in dynamic systems. 


Recently, the emergence of chaotic flows was shown to occur in the downstream pressure of a pipeline that transports highly heterogeneous mixtures of gases ~(\cite{baker2023chaos}).  The study was motivated from previous observations that 
a sudden injection of hydrogen into a natural gas pipeline could cause pressure, density, and mass flux values to diverge significantly from their nominal natural gas values ~(\cite{melaina2013blending, eames2022injection, zlotnik2023effects}).  Modeling flow phenomena for heterogeneous gas mixtures is considerably more complex than for a single gas.  The flow of natural gas in a pipe is governed by a state equation and two PDEs that describe mass and momentum conservation.  Another PDE must be added to track the concentration of each additional gas constituent.  Because societies worldwide are investing in clean energy resources for a transition to zero carbon economies, and hydrogen is a compelling energy carrier that can be produced using renewable electricity ~(\cite{gotz2016renewable, ozturk2021comprehensive}), hydrogen blending into gas pipelines is being evaluated to decarbonize energy systems while using capital investments in infrastructure for their planned lifetimes ~(\cite{hafsi2019computational,subani2017leak,chaczykowski2018gas,elaoud2017numerical}).  While model predictive optimal control methods for pipeline operations have been extended to the setting of hydrogen blending ~(\cite{baker2023optimal}), the possibility of chaotic pressure waves would challenge contemporary pipeline operation methods ~(\cite{baker2023chaos}), which use gas compressors and pressure regulators to maintain predictable flows.  This compels methods to characterize the onset of chaotic dynamics in pipeline systems, and control techniques to suppress such responses. 


In this study, we review techniques for characterizing regions in the periodic forcing parameter space that separate chaotic and non-chaotic responses to variations of hydrogen mass fraction in the gas mixture entering a pipeline.  Sinusoidal variation of hydrogen concentration at the boundary of a pipeline can create chaotic flows with large spike surges in pressure, density, mass flux, and concentration variables that propagate along the pipe, and this was demonstrated according to the rigorous mathematical definition of chaos ~(\cite{baker2023chaos}).  Here, we apply proportional feedback control to compressor settings in order to regulate the mass flow at the pipe boundary, and demonstrate with numerical simulations that this controller can suppress chaotic dynamics along the pipe.



The rest of this paper is organized as follows.  In Section \ref{sec:pipe_flow}, the PDE system that governs the flow of a mixture of gases in a pipe is introduced and reduced to a finite-dimensional ordinary differential equation (ODE) system by applying Chebyshev discretization in space.  The method that we use to quantify chaotic dynamics in the reduced ODE system is described in Section \ref{sec:quantify_chaos}.  In Section \ref{sec:control_chaos}, we define the proportional feedback control formula, in addition to mathematically defining the suppression of chaotic behavior.  Finally, Section \ref{sec:control_chaos} describes numerical experiments in which we demonstrate the performance of the controller to suppress chaos over the space of periodic forcing parameters.  

\section{Heterogeneous Gas Flow in a Pipe} \label{sec:pipe_flow}

The time-dependent flow of a mixture of natural gas and hydrogen in a horizontal pipe is modeled with the parabolic PDE system ~(\cite{chaczykowski2018gas, baker2023chaos})
\begin{eqnarray}
\frac{\partial \rho^{(m)}}{\partial t} +\frac{\partial}{\partial x}\left(\frac{\rho^{(m)}}{\rho^{(1)}+\rho^{(2)}} \varphi \right) &=&0, \label{eq:flow1} \\
\frac{\partial}{\partial x}\left(\sigma^2_1\rho^{(1)}+\sigma^2_2\rho^{(2)} \right) &=& -\frac{\lambda}{2D}\frac{ \varphi|  \varphi|}{\rho^{(1)}+\rho^{(2)}}, \label{eq:flow2} \quad
\end{eqnarray}
where (\ref{eq:flow1}) is defined for each constituent gas ($m=1$ for natural gas and $m=2$ for hydrogen).   The variables are natural gas density $\rho^{(1)}(t,x)$, hydrogen density $\rho^{(2)}(t,x)$, and total mass flux $\varphi(t,x)$ for $t\in [0,T]$ and $x\in[0,\ell]$, where $T$ denotes the time horizon and $\ell$ denotes the length of the pipe. The diameter and friction factor of the pipe are denoted by $D$ and $\lambda$. The wave speeds in natural gas and hydrogen are denoted by $\sigma_1$ and $\sigma_2$, respectively.
From Dalton's law, the total pressure $p$ of the mixture of hydrogen and natural gas is equal to the summation of their partial pressures.  For an ideal equation of state, the total pressure is given by
$ p=\sigma_1^2\rho^{(1)}+\sigma_2^2\rho^{(2)}$.  The dynamics in (\ref{eq:flow2}) assume this ideal relation.  Because the ideal equation of state qualitatively reflects the dynamics seen in the non-ideal case, we would expect to observe phenomenologically similar flow behavior with the use of more accurate nonlinear equations of state ~(\cite{botros2022performance}).  

Boundary conditions for equations (\ref{eq:flow1})-(\ref{eq:flow2}) are defined by
\begin{eqnarray}
\rho^{(m)}(t,0)&=&u(t)s^{(m)}(t), \qquad (m=1, \;2), \label{eq:bc1} \\
\varphi(t,\ell)&=& \overline q, \label{eq:bc}
\end{eqnarray}
where $s^{(m)}$ ($m=1$ and $m=2$) are known functions that represent the partial densities of the respective gases at the inlet source, and $\overline q$ is a known constant value of the total mass flux of the mixture being withdrawn from the outlet of the pipe. The function $u$ is the boundary control variable that represents the pressure increase after flow through a compressor station located at the pipe inlet.  
To more clearly distinguish between chaos and usual transient flow behavior, we specify initial conditions and all boundary parameters, except for concentration, to be as steady as possible.   Thus we choose the source pressure $\overline p=\sigma_1^2s^{(1)}+\sigma_1^2s^{(1)}$ and the withdrawal flux $\overline q$ to be constant values in time.  Moreover, the initial condition is chosen to be the steady-state solution.  As such, time-dependencies of the flow variables result only from forcing the concentration to fluctuate at the boundary.  Hydrogen concentration (mass fraction) at the boundary is defined by $\gamma=s^{(2)}/(s^{(1)}+s^{(2)})$.  If $\gamma$ and $\overline p$ are specified, then we may use their defining equations to determine necessary $s^{(m)}$ that would achieve such an outcome.  Henceforth, we assume that $\overline p$ and $\gamma$ are specified, and we define
\begin{eqnarray}
s^{(1)}(t)=(1-\gamma(t)) \overline p/ \sigma_1^{2}, \qquad
s^{(2)}(t)=\gamma(t) \overline p/ \sigma_2^{2}. \label{eq:bc_new}
\end{eqnarray} 
These definitions are enforced in the boundary conditions.

In general, the solution of (\ref{eq:flow1})-(\ref{eq:bc}) must be obtained numerically.  We apply a Chebyshev discretization scheme over the space domain $[0,\ell]$.  The interval $[0,\ell]$ is discretized at the translated extrema nodes of Chebychev polynomials defined by $x_i=\ell/2(1-\cos(i \pi/M))$ for $i=0,\dots,M$.  Let us define the state variables $\bm \rho^{(m)}_i(t)=\rho^{(m)}(t,x_i)$ and $\bm \varphi_i(t)=\varphi(t,x_i)$, along with the state vectors $\bm \rho^{(m)}=(\bm \rho^{(m)}_0,\dots,\bm \rho^{(m)}_M)^T$ and $\bm \varphi=(\bm \varphi_0,\dots,\bm \varphi_M)^T$.  
It can be shown (e.g., see \cite{ascher2011first}) that the spatial derivative of $\rho^{(m)}(t,x)$ evaluated at $x=x_i$ is approximately equal to the $(i+1)$-th column of $\bm D \bm \rho^{(m)}(t)$, where $\bm D$ is the $(M+1)\times (M+1)$ differentiation matrix defined component-wise by
\begin{eqnarray*}
        \bm D_{ij} &=& \sum_{n\not= j} \frac{1}{x_j-x_n}, \qquad (i=j), \\
\bm D_{ij} &=& \frac{1}{x_j-x_i}\prod_{ n\neq i,j} \frac{x_i-x_n}{x_j-x_n}, \qquad (i\not=j).
\end{eqnarray*}
Similar approximations are made for spatial derivatives of other flow variables.  By discretizing the PDEs in (\ref{eq:flow1})-(\ref{eq:flow2}) at the collocation nodes, we obtain the system of ODEs in the vector variables $\bm \rho^{(1)}$, $\bm \rho^{(2)}$, and $\bm \varphi$ given by
\begin{eqnarray}
 \frac{d \bm{\rho}^{(m)}}{dt} +\bm D \left(\frac{ \bm \rho^{(m)}}{\bm \rho^{(1)}+\bm \rho^{(2)}} \odot \bm \varphi \right) &=&0, \label{eq:cheb1} \\
\bm D \left(\sigma^2_1 \bm \rho^{(1)}+\sigma^2_2 \bm \rho^{(2)} \right) &=& -\frac{\lambda}{2D}\frac{ \bm \varphi \odot |  \bm \varphi|}{\bm \rho^{(1)}+\bm \rho^{(2)}}, \label{eq:cheb2} \quad
\end{eqnarray}
where $\odot$ denotes the component-wise product, and quotients are defined component-wise as well.  The boundary conditions are integrated into the ODEs with the substitutions $\bm \rho^{(m)}_0(t)=s^{(m)}(t)$ and $\bm \varphi_M(t)=\overline q$.  The solutions of (\ref{eq:cheb1})-(\ref{eq:cheb2}) are obtained numerically with Matlab using the stiff ODE solver \verb"ode15s" ~(\cite{shampine1997matlab}) using the collocation points $\{x_i\}$ corresponding to $M=16$.  A nontrivial mass matrix is included in the obvious way.
 
\section{Method to Quantify Chaos in Pipe Flow} \label{sec:quantify_chaos}

As discussed above, the periodic forcing function $\gamma(t)$ that determines partial pressures at the pipe inlet according to (\ref{eq:bc_new}) is the only time-varying  parameter.  We decompose $\gamma$ into its frequency components and analyze the effect that each individual component has on the solution.  Let us define a sinusoidal variation in hydrogen mass fraction as
\begin{equation} \label{eq:gamma_def}
    \gamma(t)= \overline \gamma \left(1+\kappa \sin(2\pi\omega t)  \right),
\end{equation}
where $\omega$ is the frequency of the sinusoid, $\kappa$ is the amplitude factor, and $\overline \gamma$ is the mean concentration of hydrogen about which the sinusoid oscillates.  Our previous study ~(\cite{baker2023chaos}) examined partitions of the boundary parameter plane $(\omega,\kappa)$ that resulted in monotonic, periodic, and chaotic response regions. Below, we design a controller to suppress chaotic responses and eliminate the area in the $(\omega,\kappa)$ plane that defines the chaotic response region.   We suppose that a control design is successful at suppressing chaos if it results in an empty chaotic response region when applied.  Before we present the control formula in the following section, we first discuss how chaos in solutions to the reduced system (\ref{eq:cheb1})-(\ref{eq:cheb2}) can be rigorously quantified.   

We consider a solution to be chaotic if it is highly sensitive to initial conditions ~(\cite{lorenz1995essence}).  The extent of chaos in a finite-dimensional system can be quantified by the largest Lyapunov exponent of the system ~(\cite{benettin1980lyapunov}).  This measure provides an estimate on the exponential rate of divergence between two solutions that begin their evolution relatively close to one another.  Numerical methods have been developed to approximate the largest Lyapunov exponent of a finite-dimensional system ~(\cite{wolf1985determining, rosenstein1993practical, brown1991computing}).  These methods require an appropriate embedding dimension, which could be problematic for solutions that are generated by an underlying infinite-dimensional PDE system.  Therefore, we measure the extent of chaos in the reduced order system by using a rate of divergence between two specific solutions.  Note that this measure does not necessarily provide an estimate of the largest Lyapunov exponent of the continuum system.  However, because the discretized system (\ref{eq:cheb1})-(\ref{eq:cheb2}) is an approximation of the continuum system (\ref{eq:flow1})-(\ref{eq:flow2}) that retains the strongly dissipative term on the right-hand side, we suppose that observations of chaos in the former system indicate chaos in the latter.

Consider two numerical solutions $\psi_1(t_n)=\psi(t_n,\ell)$ and $\psi_2(t_n)=\psi(t_n,\ell)$ with $|\psi_2(0)-\psi_1(0)|<\delta$, where $\delta>0$ is small relative to $\psi_1(0)$. Here, $\psi(t_n,\ell)$ represents the value of any one of the numerical flow variables evaluated at the pipe outlet $x=\ell$ at time $t_n=(n/N)T$ for $n=0,\dots,N$.  As in our previous study ~(\cite{baker2023chaos}), the exponential divergence between the two solutions is quantified by
\begin{eqnarray}
    \mathcal C_{\psi}  &=&\frac{1}{|I_T|}\sum_{t_n\in I_T}\log|\psi_2(t_n)-\psi_1(t_n)| \nonumber \\
                && \qquad   -  \frac{1}{|I_0|}\sum_{t_n\in I_0}\log|\psi_2(t_n)-\psi_1(t_n)|, \quad \label{eq:diverg}
\end{eqnarray}
with $I_0=[t_{n_0},t_{n_1}]$ and $I_T=[t_{n_2},t_{n_3}]$, where $n_0<n_1<n_2<n_3$.
If one flow variable shows chaotic behavior, then we would expect the other flow variables to show similar chaotic behavior. Because the measures $\mathcal C_{\rho^{(1)}}$, $\mathcal C_{\rho^{(2)}}$, and $\mathcal C_{p}$ are each sensitive to parameters and scalings, we define the \emph{chaos measure}
\begin{eqnarray}
     \mathcal C = \min(\mathcal C_{\rho^{(1)}},\mathcal C_{\rho^{(2)}},\mathcal C_{p}). \label{eq:chaos}
\end{eqnarray}
Because of the sensitivity of individual measures, we suppose that the minimum of these measures provides more sensitivity to the presence of chaos than using any one individually.  Large and positive $\mathcal C$ indicates exponential divergence between the two solutions over some time interval.  However, small and positive $\mathcal C$ does not necessarily imply that the solutions diverge from one another.  The \emph{chaotic interface} in the plane $(\omega,\kappa)$ is defined to be the function $\kappa_*(\omega)$, where for each $\omega$, $\kappa_*(\omega)$ is the lower bound that satisfies $\mathcal C(\omega,\kappa)> C$ for all $\kappa > \kappa_*(\omega)$. We find that $C=0.5$ is a useful threshold to identify chaos for the system parameters considered in this paper, and we use this threshold moving forward.  We say that the solution corresponding to the forcing parameter pair $(\omega,\kappa)$ is chaotic if $\kappa>\kappa_*(\omega)$.

\section{Control Design to Suppress Chaos} \label{sec:control_chaos}

 \begin{figure*}[ht]
    \centering
    \begin{subfigure}[t]{0.3\textwidth}
        \centering
        \includegraphics[width=1.0\linewidth]{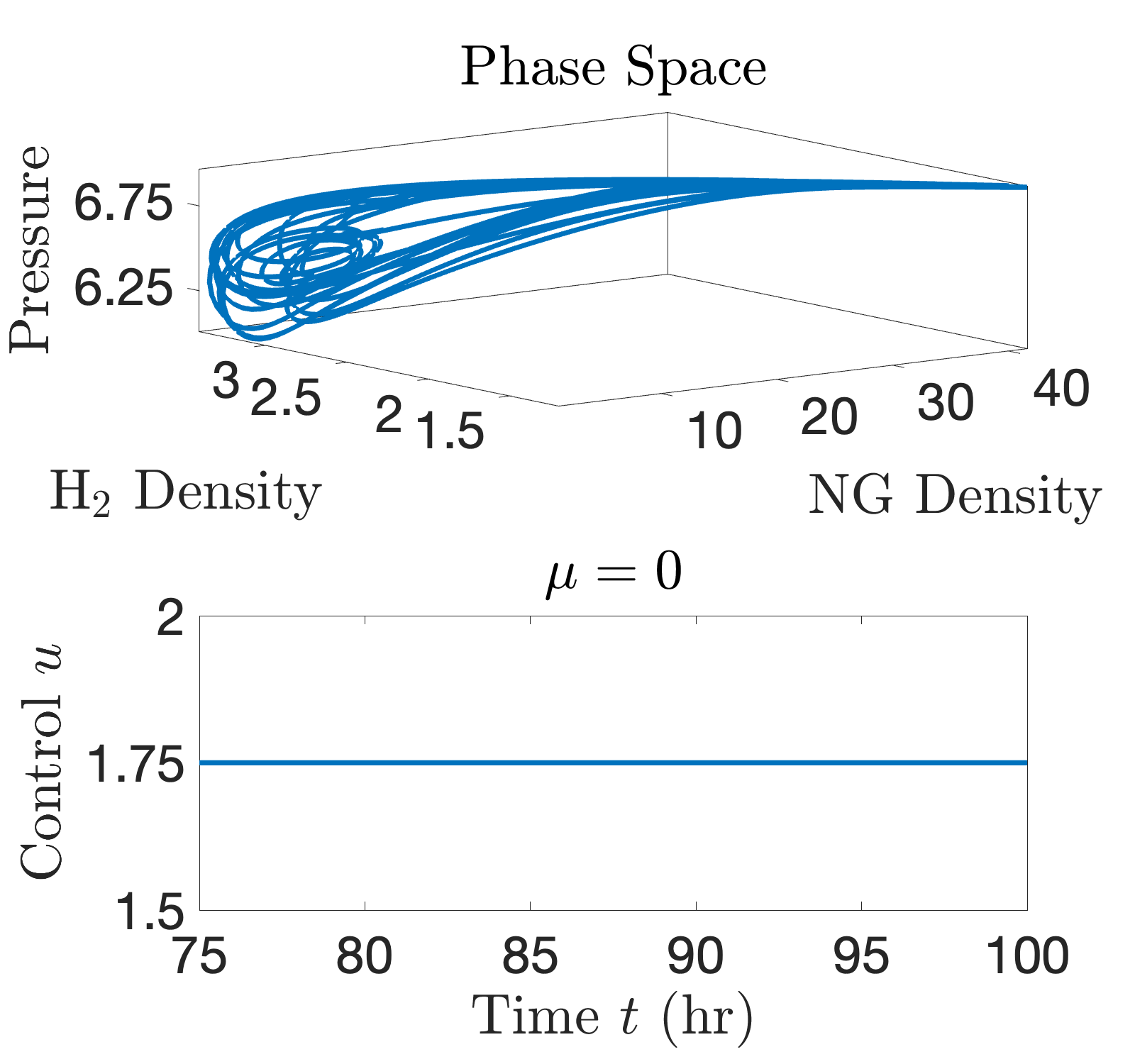}
    \end{subfigure}%
    ~ 
    \begin{subfigure}[t]{0.3\textwidth}
        \centering
        \includegraphics[width=1.0\linewidth]{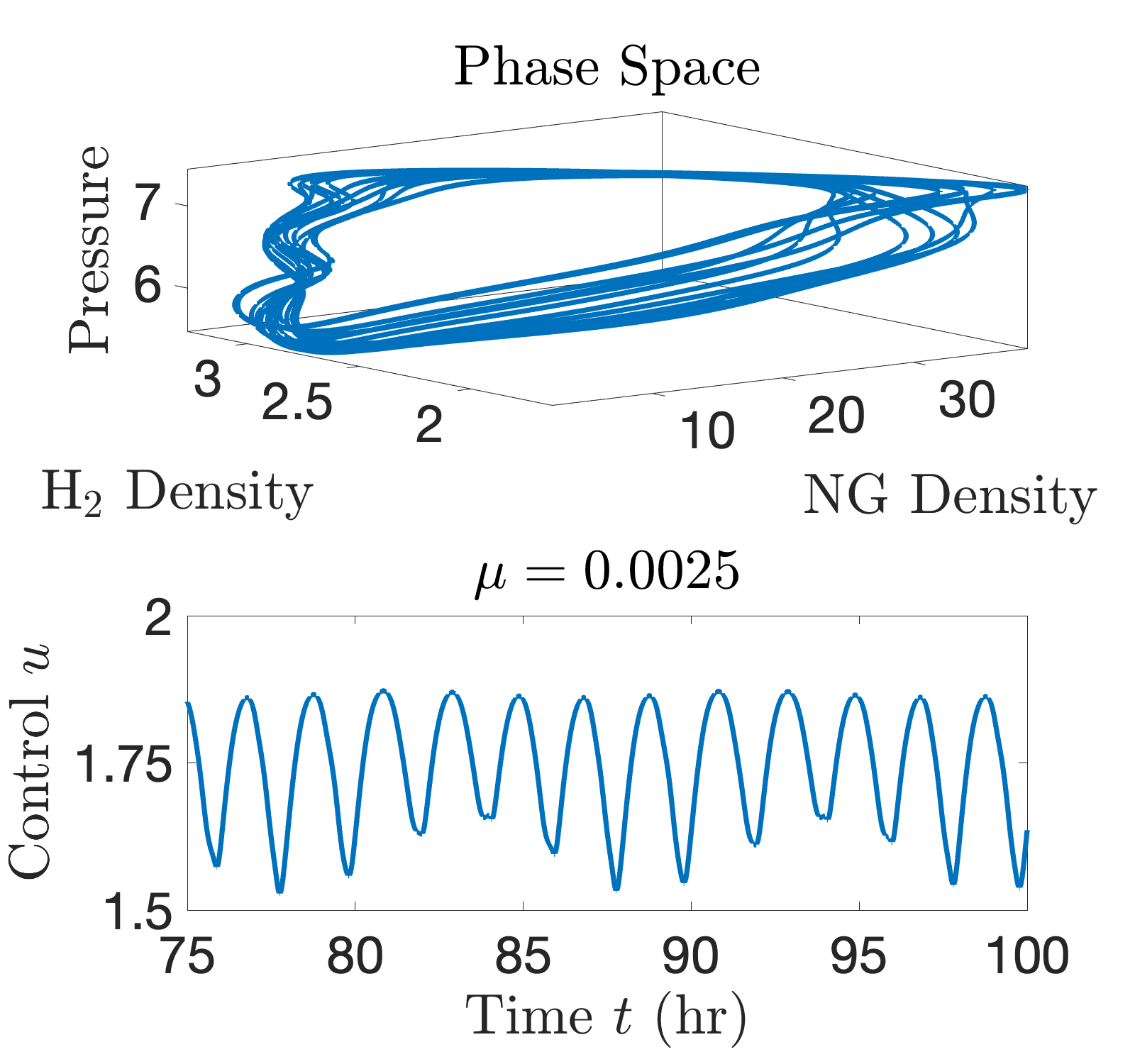}
    \end{subfigure}
        ~ 
    \begin{subfigure}[t]{0.3\textwidth}
        \centering
        \includegraphics[width=1.0\linewidth]{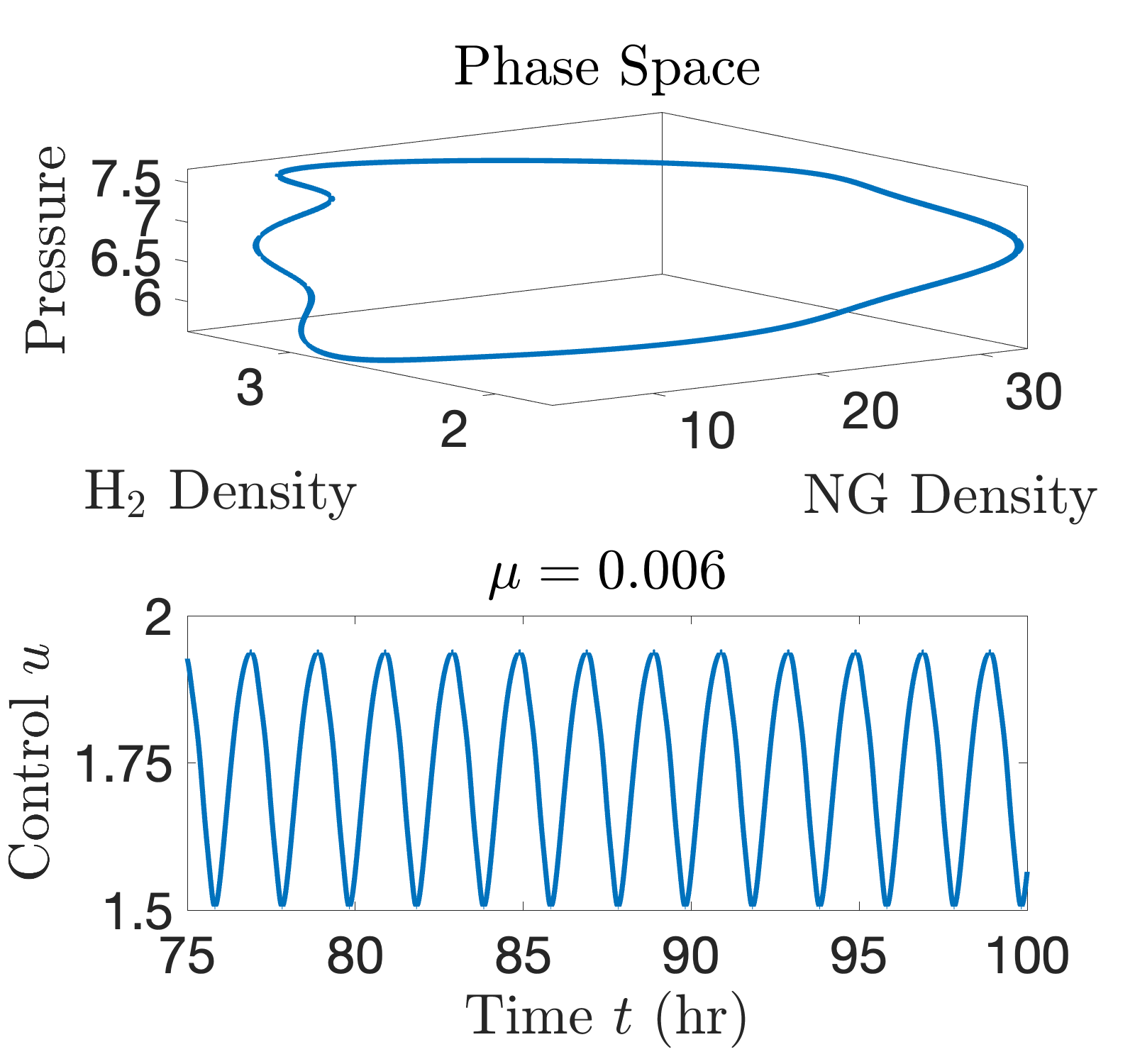}
    \end{subfigure}
    \caption{Suppressing chaos by increasing the control gain $\mu$ in (\ref{eq:control}).  Flow variables in phase space are evaluated at the pipe outlet $x=\ell$ for $t\in [75,100]$.  The concentration boundary condition parameters are $\omega=0.5$ cyc/hr, $\kappa=0.85$, and $\overline \gamma=0.2$. The remaining parameters are $\overline \mu=1.75$, $\overline p=4$ MPa, $\overline q=75$ kg/m$^2$s, $T=100$ hr, $\ell=50$ km, $D=0.5$ m, $\lambda=0.011$, $\sigma_1=338$ m/s, and $\sigma_2=4\sigma_1$.}
    \label{fig:escape_chaos}
\end{figure*}

 \begin{figure*}[ht]
    \centering
    \begin{subfigure}[t]{0.3\textwidth}
        \centering
        \includegraphics[width=1.0\linewidth]{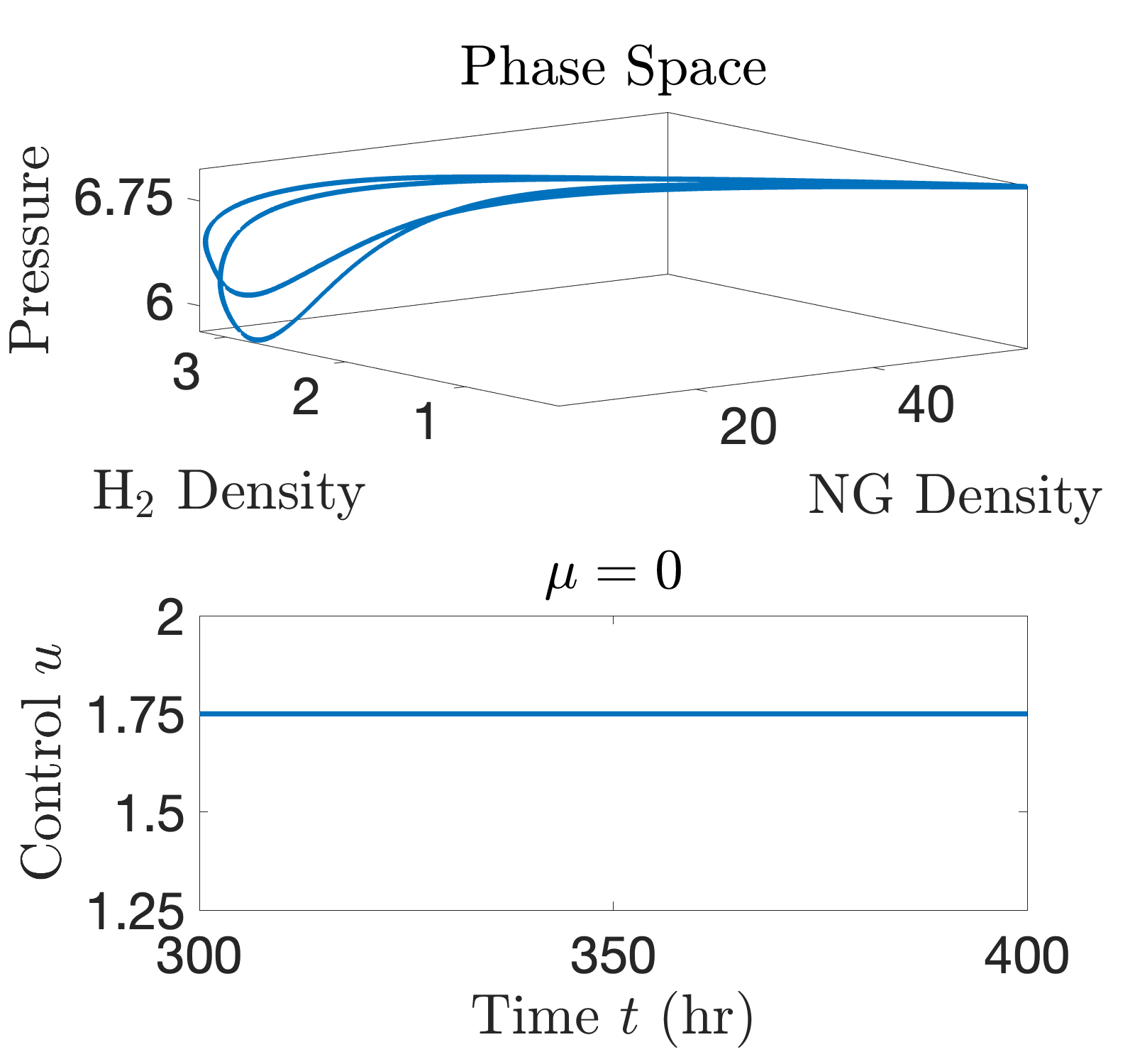}
    \end{subfigure}%
    ~ 
    \begin{subfigure}[t]{0.3\textwidth}
        \centering
        \includegraphics[width=1.0\linewidth]{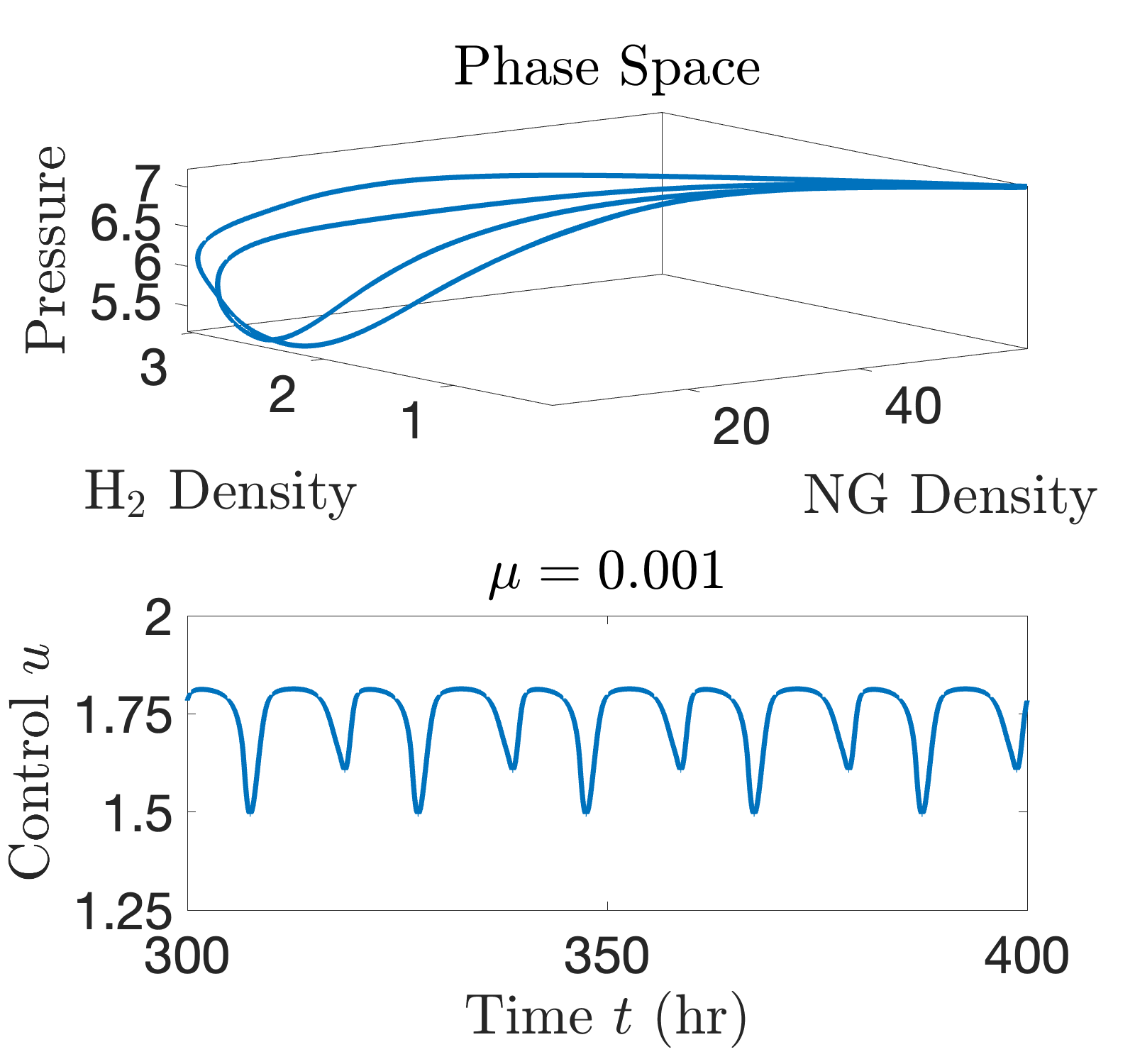}
    \end{subfigure}
        ~ 
    \begin{subfigure}[t]{0.3\textwidth}
        \centering
        \includegraphics[width=1.0\linewidth]{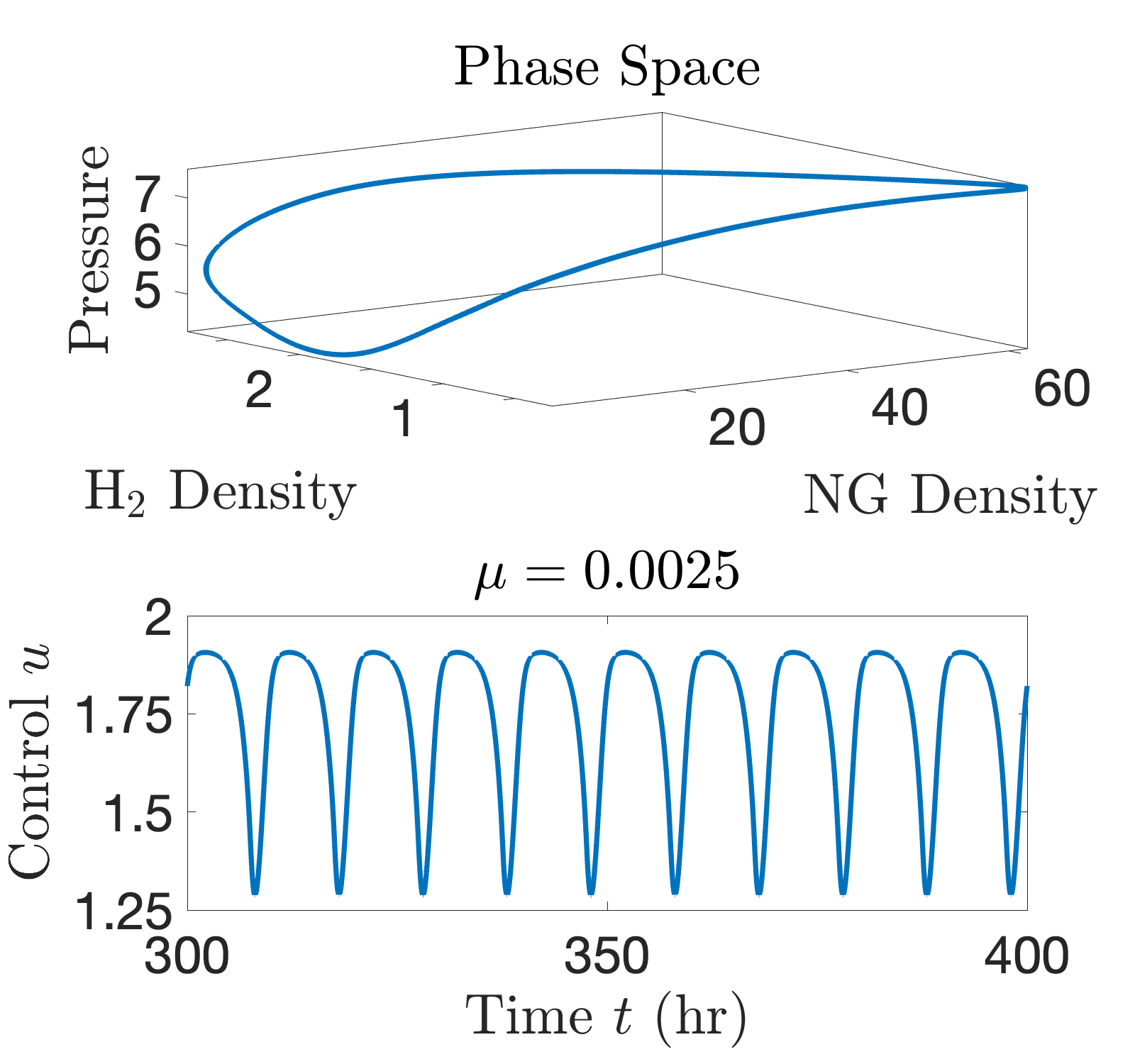}
    \end{subfigure}
    \caption{Suppressing period doubling bifurcations.  The time horizon is $T=400$ hr and the concentration parameters are $\omega=0.1$ cyc/hr, $\kappa=0.98$, and $\overline \gamma=0.2$. The remaining parameters are the same as those in Fig. \ref{fig:escape_chaos}.}
    \label{fig:escape_bifurcation}
\end{figure*}

Several different control designs varying in complexity may each demonstrate a capability of suppressing chaos.  The design that we propose is based on proportional feedback.  Our motivation for feedback control comes from industrial accessibility.  That is, because pipelines have sensors that record time series measurements of pressure and mass flux values at high frequency rates, we suppose that a chaos suppression controller could be integrated into a pipeline control system.  We propose a control action to stabilize mass flow defined by
\begin{equation}
    u(t)=\overline \mu-\mu \left( \varphi(t,0)- \varphi(0,0)\right), \label{eq:control}
\end{equation}
where $\overline \mu$ and $\mu$ are scalars that represent the baseline compressor ratio and the gain of the proportional feedback controller, respectively.  When $\mu=0$, the feedback controller is deactivated and the resulting solution is called the baseline solution.  We say that the controller with gain $\mu$ suppresses chaos at the operating point $(\omega,\kappa)$ if $\mathcal C(\omega,\kappa;\mu)< 0.5$, where $\mathcal C$ is defined as in the case without feedback control ($\mu\equiv 0$). The chaotic interface, $\kappa_*(\omega;\mu)$, corresponding to the gain $\mu$ is also defined as in the case without feedback control using equations (\ref{eq:diverg}) and (\ref{eq:chaos}).  Although the control design in (\ref{eq:control}) is conceptually simple, we will demonstrate that it is highly effective at suppressing chaotic dynamics over the entire $(\omega,\kappa)$ parameter plane.  

Let us consider an example.  We select a concentration variation parameter pair $(\omega,\kappa)=(0.5,0.85)$ that results in a chaotic solution when the controller is deactivated, i.e. $\mathcal C(\omega,\kappa;0)>0.5$.  In Fig. \ref{fig:escape_chaos}, this chaotic baseline solution along with two additional simulations that correspond to two nonzero control gains are displayed to demonstrate a suppression of chaos as the control gain increases.  All of the pipeline and boundary condition parameters, except for the control gains, are unchanged between the three simulations. The phase space portrait of the baseline solution corresponding to $\mu=0$ is shown on the left-hand side of Fig. \ref{fig:escape_chaos}.  Uninfluenced by the controller, the solution in this case regularly changes direction in phase space with sharp gradients.  This type of behavior is typical for uncontrolled solutions ~(\cite{baker2023chaos}).  In the middle of Fig. \ref{fig:escape_chaos}, we show the control action for a gain of $\mu=0.0025$ and the resulting phase portrait.  By inspection, we find that the controlled solution for this gain value exhibits less hysteresis than the baseline (uncontrolled) solution, but still does not have a closed periodic orbit in the physical phase space.  Increasing the gain further to $\mu=0.006$, we obtain a periodic phase portrait with smooth gradients and a stabilizing periodic control associated with this solution, as shown on the right-hand side of Fig. \ref{fig:escape_chaos}.  While this example demonstrates a transition to periodicity as the gain increases, this does not show quantification of chaos in the rigorous mathematical sense, i.e. that $\mathcal C(\omega,\kappa;\mu)< 0.5$ for one of the nonzero values of $\mu$.  Nevertheless, one of the results in our previous study ~(\cite{baker2023chaos}) shows that the areas in the $(\omega,\kappa)$ parameter plane that correspond to periodic and non-chaotic dynamics in fact coincide.  Therefore, we associate the transition to periodicity with the suppression of chaos in a manner that is sufficient to confirm by visual inspection of simulations. 

Although the above example may demonstrate a capability of the controller to suppress chaos and stabilize downstream flows to a periodic orbit for the operating point $(\omega,\kappa)=(0.5,0.85)$, the effectiveness of the controller does not appear to be limited to a particular operating condition.  For another example, consider the operating frequency $\omega=0.1$ cyc/hr. The baseline solution associated with this frequency approaches a periodic orbit, as seen on the left-hand side of Fig. \ref{fig:escape_bifurcation}; however, this baseline solution experiences a sequence of period doubling bifurcations as the amplitude factor $\kappa$ increases to unity.  For the amplitude factor $\kappa=0.98$, the baseline solution has twice the period of the forcing concentration $\gamma(t)$ as a result of the two loops that comprise the orbit, as seen on the left-hand side of Fig. \ref{fig:escape_bifurcation}.  Inspecting Fig. \ref{fig:escape_bifurcation}, we observe that the controller suppresses period doubling, in the sense that the two loops in the orbit of the baseline solution transition to a single loop for a sufficiently large control gain.  The examples in Figs. \ref{fig:escape_chaos}-\ref{fig:escape_bifurcation} demonstrate that the controller may have the ability to suppress chaos while simultaneously stabilizing its dynamics to a simple periodic orbit with the use of sufficiently large control gains that depend on the operating points (and other pipeline parameters).  Moreover, the examples illustrate that a stabilizing controller has the potential to significantly reduce large and undesirable gradients in phase space that could otherwise appear in the uncontrolled solution.

\begin{figure}[t!]
\centering
\includegraphics[width=1.0\linewidth]{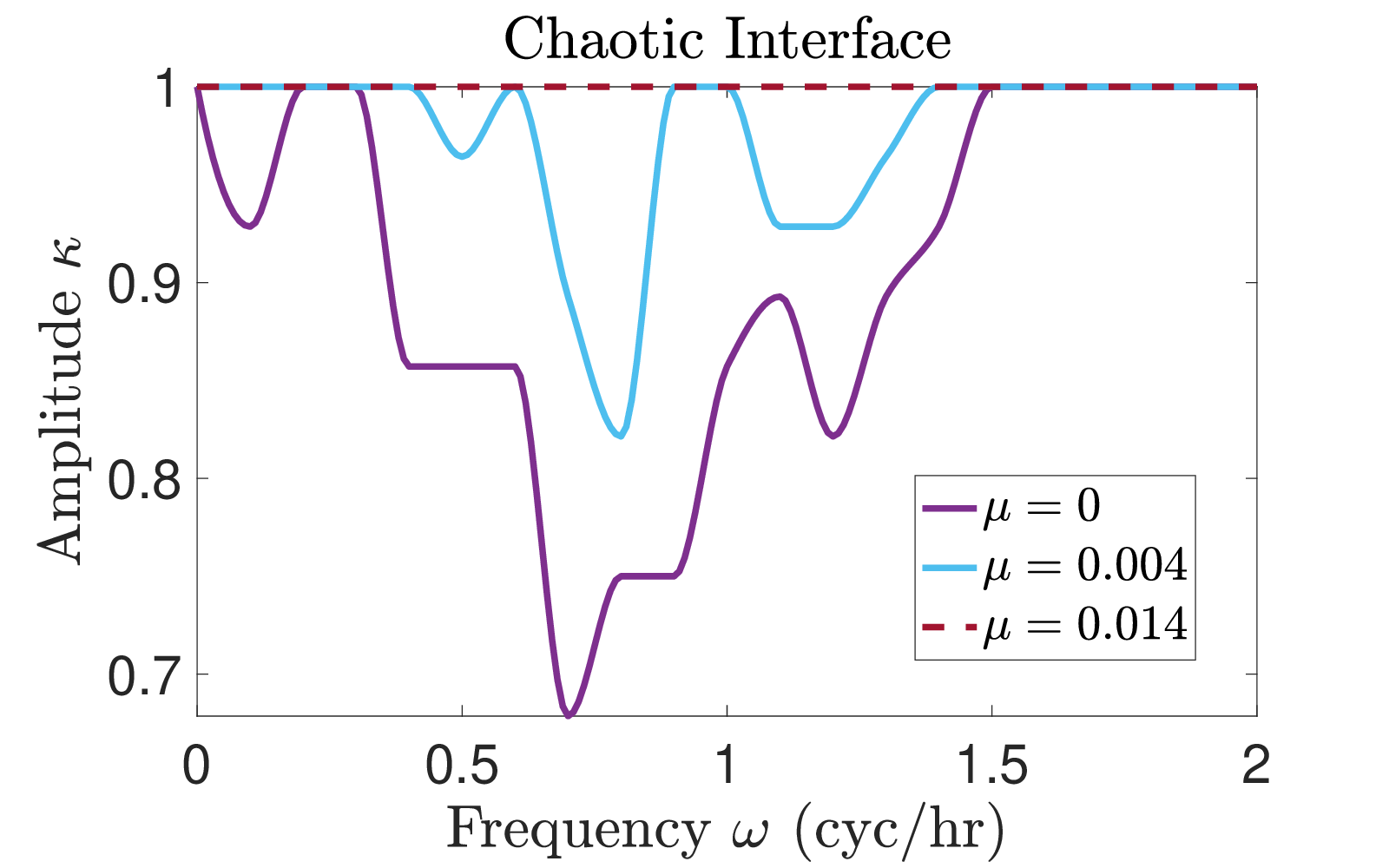}
\caption{Chaotic interfaces $\kappa_*(\omega;\mu)$ in the $(\omega,\kappa)$ plane as functions of $\omega$ for various control gains $\mu$. The region defined by $0\le \omega\le 2$ and $0.5\le \kappa \le 1$ is discretized into a $21\times 15$ grid of pairs of points.  For each pair and for each control gain, two solutions are simulated for 100 hours whose initial conditions correspond to withdrawal rates $\overline q=75$ and $\overline q=75.1$.   The value $\mathcal C(\omega,\kappa;\mu)$ is determined from the two solutions using (\ref{eq:chaos}) with intervals $I_0=[0.08N,0.15N]$ and $I_T=[0.5N,0.8N]$, where $N=10,000$.   The chaotic interface curve $\kappa_*(\omega;\mu)$ is obtained from the set of values $\mathcal C(\omega,\kappa;\mu)$ on the grid according to its definition and then interpolated with shape-preserving splines.  The remaining parameters are the same as those in Fig. \ref{fig:escape_chaos}.}
\label{fig:CI}
\end{figure}

We conclude by analyzing the chaotic interface $\kappa_*(\omega;\mu)$ as a function of frequency for various control gains.  
Figure \ref{fig:CI} depicts three chaotic interfaces in the $(\omega,\kappa)$ parameter space that correspond to the three control gains listed in the figure.   It may appear from Fig. \ref{fig:CI} that for any two chaotic interfaces, $\kappa_*(\omega;\mu_1)$ and $\kappa_*(\omega;\mu_2)$, the inequality $\kappa_*(\omega;\mu_1)\le \kappa_*(\omega;\mu_2)$ is satisfied for all $\omega$ whenever $\mu_1\le \mu_2$.  However, from extensive numerical simulations, we note that there are cases for which this proposed inequality does not hold.  Such non-monotonicity in computed criteria may be a consequence of sensitivity issues that arise in characterizing chaotic dynamics numerically.  Nonetheless, the chaotic interfaces shown in Fig. \ref{fig:CI} are in agreement with the conclusions from the examples presented above.  The interfaces successfully demonstrate that a sufficiently large control gain $\mu$ in a proportional feedback controller of the form in (\ref{eq:control})  could in fact be capable of suppressing chaos over the entire ensemble of operating points that comprise the boundary parameter plane.  

\section{Conclusions} \label{sec:conclusion}

We have reviewed a new problem in the boundary control of chaotic dynamics in a distributed-parameter system defined by partial differential equations.  This setting has a compelling application to suppressing the chaotic flows that can potentially occur in the downstream pressure of a natural gas pipeline into which a time-varying mass fraction of hydrogen gas is injected.  We demonstrate a simple single-pipe test system in which the undesirable chaotic dynamics emerge, and then propose a local controller that modulates the pipe inlet pressure in order to stabilize the inlet flow using proportional feedback.  This controller is shown to completely suppress chaotic dynamics in flows and pressures throughout the pipe for periodic hydrogen mass fraction variations in the entire phase plane of frequency and amplitude parameters.  Existing compressor station control systems already incorporate the capability to stabilize through-flow rates, in addition to a variety of pressure control settings.  Applying the type of controller proposed in our study, natural gas pipelines could accommodate injections of hydrogen produced using renewable energy to facilitate decarbonization while using capital investments in infrastructure for their planned lifetimes.


\bibliography{ifacconf}             

\begin{thebibliography}{32}
\providecommand{\natexlab}[1]{#1}
\providecommand{\url}[1]{\texttt{#1}}
\providecommand{\urlprefix}{URL }
\expandafter\ifx\csname urlstyle\endcsname\relax
  \providecommand{\doi}[1]{doi:\discretionary{}{}{}#1}\else
  \providecommand{\doi}{doi:\discretionary{}{}{}\begingroup \urlstyle{rm}\Url}\fi

\bibitem[{Alipoor and Mazaheri(2016)}]{alipoor2016combustion}
Alipoor, A. and Mazaheri, K. (2016).
\newblock Combustion characteristics and flame bifurcation in repetitive extinction-ignition dynamics for premixed hydrogen-air combustion in a heated micro channel.
\newblock \emph{Energy}, 109, 650--663.

\bibitem[{Ascher and Greif(2011)}]{ascher2011first}
Ascher, U.M. and Greif, C. (2011).
\newblock \emph{A first course on numerical methods}.
\newblock SIAM.

\bibitem[{Baker et~al.(2023{\natexlab{a}})Baker, Kazi, Platte, and Zlotnik}]{baker2023optimal}
Baker, L.S., Kazi, S.R., Platte, R.B., and Zlotnik, A. (2023{\natexlab{a}}).
\newblock Optimal control of transient flows in pipeline networks with heterogeneous mixtures of hydrogen and natural gas.
\newblock In \emph{2023 American Control Conference (ACC)}, 1221--1228. IEEE.

\bibitem[{Baker et~al.(2023{\natexlab{b}})Baker, Kazi, and Zlotnik}]{baker2023chaos}
Baker, L.S., Kazi, S.R., and Zlotnik, A. (2023{\natexlab{b}}).
\newblock Transitions from monotonicity to chaos in gas mixture dynamics in pipeline networks.
\newblock \emph{PRX Energy}, 2, 033008.

\bibitem[{Benettin et~al.(1980)Benettin, Galgani, Giorgilli, and Strelcyn}]{benettin1980lyapunov}
Benettin, G., Galgani, L., Giorgilli, A., and Strelcyn, J.M. (1980).
\newblock Lyapunov characteristic exponents for smooth dynamical systems and for hamiltonian systems; a method for computing all of them. part 1: Theory.
\newblock \emph{Meccanica}, 15, 9--20.

\bibitem[{Boccaletti et~al.(2000)Boccaletti, Grebogi, Lai, Mancini, and Maza}]{boccaletti2000control}
Boccaletti, S., Grebogi, C., Lai, Y.C., Mancini, H., and Maza, D. (2000).
\newblock The control of chaos: theory and applications.
\newblock \emph{Physics reports}, 329(3), 103--197.

\bibitem[{Botros and Jensen(2022)}]{botros2022performance}
Botros, K.K. and Jensen, L. (2022).
\newblock Performance of twelve different equations of state for natural gas and hydrogen blends.
\newblock In \emph{International Pipeline Conference}, volume 86564, V001T08A002. American Society of Mechanical Engineers.

\bibitem[{Brown et~al.(1991)Brown, Bryant, and Abarbanel}]{brown1991computing}
Brown, R., Bryant, P., and Abarbanel, H.D.I. (1991).
\newblock Computing the lyapunov spectrum of a dynamical system from an observed time series.
\newblock \emph{Physical review A}, 43(6), 2787.

\bibitem[{Chaczykowski et~al.(2018)Chaczykowski, Sund, Zarodkiewicz, and Hope}]{chaczykowski2018gas}
Chaczykowski, M., Sund, F., Zarodkiewicz, P., and Hope, S.M. (2018).
\newblock Gas composition tracking in transient pipeline flow.
\newblock \emph{Journal of Natural Gas Science and Engineering}, 55, 321--330.

\bibitem[{Ding et~al.(2017)Ding, Song, Yang, Litak, Wang, Yao, and Ma}]{ding2017analysis}
Ding, S.L., Song, E.Z., Yang, L.P., Litak, G., Wang, Y.Y., Yao, C., and Ma, X.Z. (2017).
\newblock Analysis of chaos in the combustion process of premixed natural gas engine.
\newblock \emph{Applied Thermal Engineering}, 121, 768--778.

\bibitem[{Dudnik et~al.(1983)Dudnik, Kuznetsov, Minakova, and Romanovskii}]{dudnik1983synchronization}
Dudnik, E.N., Kuznetsov, Y.I., Minakova, I., and Romanovskii, Y.M. (1983).
\newblock Synchronization in systems with strange attractors.
\newblock \emph{Moscow Univ. Phys. Bull. Ser}, 3(24), 84--87.

\bibitem[{Eames et~al.(2022)Eames, Austin, and Wojcik}]{eames2022injection}
Eames, I., Austin, M., and Wojcik, A. (2022).
\newblock Injection of gaseous hydrogen into a natural gas pipeline.
\newblock \emph{International Journal of Hydrogen Energy}, 47(61), 25745--25754.

\bibitem[{Eisenman et~al.(2005)Eisenman, Yu, and Tziperman}]{eisenman2005westerly}
Eisenman, I., Yu, L., and Tziperman, E. (2005).
\newblock Westerly wind bursts: Enso’s tail rather than the dog?
\newblock \emph{Journal of Climate}, 18(24), 5224--5238.

\bibitem[{Elaoud et~al.(2017)Elaoud, Hafsi, and Hadj-Taieb}]{elaoud2017numerical}
Elaoud, S., Hafsi, Z., and Hadj-Taieb, L. (2017).
\newblock Numerical modelling of hydrogen-natural gas mixtures flows in looped networks.
\newblock \emph{Journal of Petroleum Science and Engineering}, 159, 532--541.

\bibitem[{G{\"o}tz et~al.(2016)G{\"o}tz, Lefebvre, M{\"o}rs, Koch, Graf, Bajohr, Reimert, and Kolb}]{gotz2016renewable}
G{\"o}tz, M., Lefebvre, J., M{\"o}rs, F., Koch, A.M., Graf, F., Bajohr, S., Reimert, R., and Kolb, T. (2016).
\newblock Renewable power-to-gas: A technological and economic review.
\newblock \emph{Renewable energy}, 85, 1371--1390.

\bibitem[{Hafsi et~al.(2019)Hafsi, Elaoud, and Mishra}]{hafsi2019computational}
Hafsi, Z., Elaoud, S., and Mishra, M. (2019).
\newblock A computational modelling of natural gas flow in looped network: Effect of upstream hydrogen injection on the structural integrity of gas pipelines.
\newblock \emph{Journal of Natural Gas Science and Engineering}, 64, 107--117.

\bibitem[{Hu et~al.(1995)Hu, Qu, and He}]{hu1995feedback}
Hu, G., Qu, Z., and He, K. (1995).
\newblock Feedback control of chaos in spatiotemporal systems.
\newblock \emph{International Journal of Bifurcation and Chaos}, 5(04), 901--936.

\bibitem[{Kocarev(2001)}]{kocarev2001chaos}
Kocarev, L. (2001).
\newblock Chaos-based cryptography: a brief overview.
\newblock \emph{IEEE Circuits and Systems Magazine}, 1(3), 6--21.

\bibitem[{Lorenz(1995)}]{lorenz1995essence}
Lorenz, E.N. (1995).
\newblock \emph{The essence of chaos}.
\newblock University of Washington press.

\bibitem[{Melaina et~al.(2013)Melaina, Antonia, and Penev}]{melaina2013blending}
Melaina, M.W., Antonia, O., and Penev, M. (2013).
\newblock Blending hydrogen into natural gas pipeline networks: A review of key issues.
\newblock Technical report.

\bibitem[{Ott et~al.(1990)Ott, Grebogi, and Yorke}]{ott1990controlling}
Ott, E., Grebogi, C., and Yorke, J.A. (1990).
\newblock Controlling chaos.
\newblock \emph{Physical review letters}, 64(11), 1196.

\bibitem[{Ozturk and Dincer(2021)}]{ozturk2021comprehensive}
Ozturk, M. and Dincer, I. (2021).
\newblock A comprehensive review on power-to-gas with hydrogen options for cleaner applications.
\newblock \emph{International Journal of Hydrogen Energy}, 46(62), 31511--31522.

\bibitem[{Pizza et~al.(2008)Pizza, Frouzakis, Mantzaras, Tomboulides, and Boulouchos}]{pizza2008dynamics}
Pizza, G., Frouzakis, C.E., Mantzaras, J., Tomboulides, A.G., and Boulouchos, K. (2008).
\newblock Dynamics of premixed hydrogen/air flames in microchannels.
\newblock \emph{Combustion and Flame}, 152(3), 433--450.

\bibitem[{Pyragas(1992)}]{pyragas1992continuous}
Pyragas, K. (1992).
\newblock Continuous control of chaos by self-controlling feedback.
\newblock \emph{Physics letters A}, 170(6), 421--428.

\bibitem[{Rosenstein et~al.(1993)Rosenstein, Collins, and De~Luca}]{rosenstein1993practical}
Rosenstein, M.T., Collins, J.J., and De~Luca, C.J. (1993).
\newblock A practical method for calculating largest lyapunov exponents from small data sets.
\newblock \emph{Physica D: Nonlinear Phenomena}, 65(1-2), 117--134.

\bibitem[{Schuster and Stemmler(1997)}]{schuster1997control}
Schuster, H.G. and Stemmler, M.B. (1997).
\newblock Control of chaos by oscillating feedback.
\newblock \emph{Physical Review E}, 56(6), 6410.

\bibitem[{Shampine and Reichelt(1997)}]{shampine1997matlab}
Shampine, L.F. and Reichelt, M.W. (1997).
\newblock The matlab ode suite.
\newblock \emph{SIAM journal on scientific computing}, 18(1), 1--22.

\bibitem[{Subani et~al.(2017)Subani, Amin, and Agaie}]{subani2017leak}
Subani, N., Amin, N., and Agaie, B.G. (2017).
\newblock Leak detection of non-isothermal transient flow of hydrogen-natural gas mixture.
\newblock \emph{Journal of Loss Prevention in the Process Industries}, 48, 244--253.

\bibitem[{Tziperman et~al.(1994)Tziperman, Stone, Cane, and Jarosh}]{tziperman1994nino}
Tziperman, E., Stone, L., Cane, M.A., and Jarosh, H. (1994).
\newblock {El Ni{\~n}o} chaos: Overlapping of resonances between the seasonal cycle and the pacific ocean-atmosphere oscillator.
\newblock \emph{Science}, 264(5155), 72--74.

\bibitem[{Wolf et~al.(1985)Wolf, Swift, Swinney, and Vastano}]{wolf1985determining}
Wolf, A., Swift, J.B., Swinney, H.L., and Vastano, J.A. (1985).
\newblock Determining lyapunov exponents from a time series.
\newblock \emph{Physica D: nonlinear phenomena}, 16(3), 285--317.

\bibitem[{Zhang et~al.(2009)Zhang, Liu, and Wang}]{zhang2009controlling}
Zhang, H., Liu, D., and Wang, Z. (2009).
\newblock \emph{Controlling chaos: suppression, synchronization and chaotification}.
\newblock Springer Science \& Business Media.

\bibitem[{Zlotnik et~al.(2023)Zlotnik, Kazi, Sundar, Gyrya, Baker, Sodwatana, and Brodskyi}]{zlotnik2023effects}
Zlotnik, A., Kazi, S.R., Sundar, K., Gyrya, V., Baker, L., Sodwatana, M., and Brodskyi, Y. (2023).
\newblock Effects of hydrogen blending on natural gas pipeline transients, capacity, and economics.
\newblock In \emph{PSIG Annual Meeting}, PSIG--2312. PSIG.

\end{thebibliography}

\end{document}